\documentclass{amsart}
\usepackage{latexsym, amssymb, amsfonts, amsmath}

\newcommand{\cx}{\mathrm{C}(X)}
\newcommand{\cy}{\mathrm{C}(Y)}
\newcommand{\ml}{\mathbf{\{}}
\newcommand{\mr}{\mathbf{\}}}
\newcommand{\e}{\epsilon}

\newcommand{\la}{\lambda}

\newtheorem{thrm}{Theorem}[section]
\newtheorem{lemma}[thrm]{Lemma}
\newtheorem{cor}[thrm]{Corollary}
\newtheorem{prop}[thrm]{Proposition}

\theoremstyle{definition}

\newtheorem{defn}[thrm]{Definition}

\begin{document}

\title[C$^*$-algebras with stable rank one]{A class of simple C$^*$-algebras with stable rank one}

\author{George A. Elliott}
\address{Department of Mathematics, University of Toronto, Toronto, Ontario, Canada M5S 2E4}
\email{elliott@math.toronto.edu}
\author{Toan M. Ho}
\address{Department of Mathematics and Statistics, York University, 4700 Keele St., Toronto, Ontario, Canada M3J 1P3}
\email{toan@mathstat.yorku.ca}
\author{Andrew S. Toms}
\address{Department of Mathematics and Statistics, York University, 4700 Keele St., Toronto, Ontario, Canada M3J 1P3}
\email{atoms@mathstat.yorku.ca}

\thanks{The authors were supported by NSERC.}

\subjclass{Primary 60F99; secondary 60H30, 60G07, 42A16.}

\keywords{AH algebras, stable rank}

\begin{abstract}
We study the limits of inductive sequences $(A_i,\phi_i)$ where each 
$A_{i}$ is a direct sum of full matrix algebras over compact metric
spaces and each partial map of $\phi_i$ is diagonal.  We give a new
characterisation of simplicity for such algebras, and apply it to
prove that the said algebras have stable rank one whenever they are
simple and unital.  Significantly, our results do not require any 
dimension growth assumption.
\end{abstract}

\maketitle

\section{Introduction}

Let $X$ and $Y$ be compact Hausdorff spaces.
A $*$-homomorphism 
\[
\phi:\mathrm{M}_{m}(\cx) \to \mathrm{M}_{nm}(\cy)
\]
is called {\it diagonal} 
if there are $n$ continuous maps $\la_i:Y \to X$ such that
\begin{displaymath}
\phi(f) =
  \left( \begin{array}{cccc}
      f \circ \la_{1} & 0  &  \ldots & 0   \\
      0 &  f \circ \la_{2} & \ldots & 0  \\
      \vdots & \vdots &  \ddots &  \vdots \\
      0  &   0 & \ldots &      f \circ \la_{n}
       \end{array} \right).
\end{displaymath}
The $\la_i$ are called the {\it eigenvalue maps} or simply {\it eigenvalues} of $\phi$.
The multiset $\ml \la_{1},{\la}_{2}, \ldots ,{\la}_{n} \mr $is called the 
{\it eigenvalue pattern} of $\phi$ and is denoted by $\textrm{ep}(\phi)$. 
This definition can be extended to $*$-homomorphisms 
\[
\phi: \bigoplus_{i=1}^n \mathrm{M}_{n_i}(\mathrm{C}(X_i)) \to 
\bigoplus_{j=1}^m \mathrm{M}_{m_j}(\mathrm{C}(Y_j))
\]
by requiring, roughly, that each partial map
\[
\phi^{ij}: \mathrm{M}_{n_i}(\mathrm{C}(X_i)) \to \mathrm{M}_{m_j}(\mathrm{C}(Y_j))
\]
induced by $\phi$ be diagonal. (A precise definition can be found in Section 2.)  

C$^*$-algebras obtained as limits of inductive systems $(A_i,\phi_i)$ where 
\[
A_i = \bigoplus_{j=1}^{n_i} \mathrm{M}_{n_{i,j}}(\mathrm{C}(X_{i,j}))
\]
and each $\phi_i$ is diagonal form a rich class.  They include AF algebras, simple unital
A$\mathbb{T}$ algebras (and hence the irrational rotation algebras) (\cite{EE}), Goodearl algebras (\cite{goodearl}),
and some interesting examples of Villadsen and the third named author connected to
Elliott's program to classify amenable C*-algebras via K-theory (\cite{jes1}, \cite{To2}).  The structure
of these algebras is only well understood when they satisfy some additional
conditions such as (very) slow dimension growth or the combination of real rank zero, stable
rank one, and weak unperforation of the $\mathrm{K}_0$-group --- situations in which the strong form of
Elliott's classification conjecture can be verified (\cite{D}, \cite{Ell}, \cite{EG}, \cite{EGL}, \cite{G}).  

In this paper we give a new characterisation of simplicity for AH algebras with diagonal connecting maps. 
As a consequence we are able to prove that such algebras have stable rank one
whenever they are unital and simple.
The significance of our result derives from the fact that we make no assumptions
on the dimension growth of the algebras;  we obtain a general theorem on
the structure of algebras heretofore considered ``wild''.  As suggested by M. R\o rdam 
in his recent ICM address, it is high time we became friends with such algebras, as
opposed to treating them simply as a source of pathological examples.    

\vspace{2mm}
\noindent
{\it Acknowledgement.}  The third named author thanks Hanfeng Li, Cristian
Ivanescu, and Zhuang Niu for several helpful discussions and comments.

\section{Preliminaries}

\subsection{Basic notation}
We use $\mathrm{M}_n$ to denote the set of $n \times n$ complex matrices.  Given
a closed subset $E$ of a compact metric space $(X,d)$ and $\delta > 0$ we set
\[
B_{\delta}(E) = \{x \in X \ | \ d(E,x) < \delta \},
\]
and make the convention that $B_{\delta}(\emptyset) = \emptyset$.

\subsection{AH systems with diagonal maps}

\begin{defn}\label{gendiagmap}
We will say that a unital $*$-homomorphsim
\[
\phi: \bigoplus_{i=1}^n M_{n_i}(\mathrm{C}(X_i)) \to M_k(\cy) \cong M_k \otimes \cy
\]
is diagonal if there exist natural numbers $k_1,\ldots,k_n$ such
that $\sum_i k_i = k$ and $n_i | k_i$, an embedding 
\[
\iota: \bigoplus_{i=1}^n M_{k_i} \hookrightarrow M_k,
\]
and diagonal maps
\[
\phi_i:M_{n_i}(\mathrm{C}(X_i)) \to M_{k_i} \otimes \cy
\]
such that
\[
\phi = \bigoplus_{i=1}^n \phi_i.
\]
(Notice that $k_i = 0$ is allowed.)  We will say that a unital $*$-homomorphism
\[
\phi: \bigoplus_{i=1}^n \mathrm{M}_{n_i}(\mathrm{C}(X_i)) \to 
\bigoplus_{j=1}^m \mathrm{M}_{m_j}(\mathrm{C}(Y_j))
\]
is diagonal if each restriction
\[
\phi_j: \bigoplus_{i=1}^n \mathrm{M}_{n_i}(\mathrm{C}(X_i)) \to \mathrm{M}_{m_j}(\mathrm{C}(Y_j))
\]
is diagonal.
\end{defn}

Let $A$ be the limit of the inductive sequence (${A}_{i}, \phi_{i}$), where
\begin{equation} \label{equa1.1}
A_i = \bigoplus_{t=1}^{k_i} \mathrm{M}_{n_{i,t}}(\mathrm{C}(X_{i,t})),
\end{equation}
$X_{i,t}$ is a connected compact metric space, and $n_{i,t}$ and $k_{i}$ are 
natural numbers.  Define
\[
A_{i,t} := \mathrm{M}_{n_{i,t}}(\mathrm{C}(X_{i,t})),
\]
\[
X_i := X_{i,1} \sqcup X_{i,2} \sqcup \cdots \sqcup X_{i,k_i},
\]
and
\[
\phi_{i,j} := \phi_{j-1} \circ \cdots \circ \phi_i.
\]
Let 
\[
\phi_{i,j}^{t,l}:\mathrm{M}_{n_{i,t}}(\mathrm{C}(X_{i,t})) \to \mathrm{M}_{n_{j,l}}(\mathrm{C}(X_{j,l}))
\]
and
\[
\phi_{i,j}^l:  \bigoplus_{t=1}^{k_{i}}
\mathrm{M}_{n_{i,t}}(\mathrm{C}(X_{i,t})) \to
\mathrm{M}_{n_{j,l}}(\mathrm{C}(Y_{j,l}))
\]
denote the appropriate restrictions of $\phi$.  If each $\phi_i$ is unital and diagonal, then we refer to 
$(A_i,\phi_i)$ as an {\it AH system with diagonal maps}.  The limit algebra $A$
will be called a {\it diagonal AH algebra}, and we will refer to $(A_i,\phi_i)$
as a {\it decomposition} of $A$. 

Assume that $A$ as above is diagonal.  We will view $\phi_{i,j}^{t,l}$ as a diagonal
map from  $\mathrm{M}_{n_{i,t}}(\mathrm{C}(X_{i,t}))$ into the cut-down of
$\mathrm{M}_{n_{j,l}}(\mathrm{C}(X_{j,l}))$ by $\phi_{i,j}^{t,l}(1)$.
For fixed $i$ and $j$, we will denote by $\mathrm{ep}_{ij}$ the multiset which
is the union, counting multiplicity, of the eigenvalue patterns of each $\phi_{i,j}^{t,l}$;  
$\mathrm{ep}_{ij}$ is the {\it eigenvalue pattern} of $\phi_{i,j}$;
an element of $\mathrm{ep}_{ij}$ is an {\it eigenvalue map} of $\phi_{i,j}$.
For fixed $i$, $j$, and $l$, we will denote by $\mathrm{ep}_{ij}^l$ the multiset which
is the union of the eigenvalue patters of each $\phi_{i,j}^{t,l}$.

Let us now show that the bonding maps $\phi_i$ may be assumed to be injective.
Let $(A_i,\phi_i)$ be a decomposition for a diagonal AH algebra $A$ as above.
Fix $i \in \mathbb{N}$ and $1 \leq t \leq k_i$.  For each $j > i$ and $1 \leq l \leq k_j$,
Let $X_{i,t}^{j,l}$ denote the closed subset of $X_{i,t}$ which is the union
of the images of the eigenvalue maps of $\phi_{i,j}^{t,l}$.  Put
$X_{i,t}^j = \cup_l X_{i,t}^{j,l}$, and $\tilde{X}_{i,t} = \cap_j X_{i,t}^j$.
Since $X_{i,t}^j \supseteq X_{i,t}^{j+1}$, we have that $\tilde{X}_{i,t}$ is
closed subset of $X_{i,t}$.  Define
\[
\tilde{A}_{i,t} = \mathrm{M}_{n_{i,t}}(\mathrm{C}(\tilde{X}_{i,t}))      
\]
and 
\[
\tilde{A}_i = \bigoplus_{t=1}^{k_i} \tilde{A}_{i,t}.
\]
Define diagonal maps $\tilde{\phi}_{i,i+1}^{t,l}:\tilde{A}_{i,t} \to \tilde{A}_{i+1,l}$
by replacing the eigenvalue maps of $\phi_{i,i+1}^{t,l}$ with their
restrictions to $\tilde{X}_{i+1,l}$.  Define $\tilde{\phi}_i:\tilde{A}_i \to \tilde{A}_{i+1}$
in a manner analogous to the definition of $\phi_i$.  It follows that
$(\tilde{A}_i,\tilde{\phi}_i)$ is a diagonal AH system with limit $A$,
and $\tilde{\phi}_i$ is injective by construction.  We assume from here
on that all bonding maps in diagonal AH systems are injective.

One way to construct a simple diagonal AH algebra is to ensure that for 
each $i \in \mathbb{N}$ and $x$ in a specified dense subset of $X_i$
there is some $j \geq i$ such that for each $l \in \{1,\ldots,k_j\}$ the 
diagonal map $\phi_{i,j}^{t,l}$ contains the eigenvalue map
$ev_x:X_{j,l} \to X_{i,t}$ given by $ev_x(y) = x$.  The next definition gives 
and approximate version of this situation.

\begin{defn}\label{defn1}
Say that a diagonal AH algebra $A$ with decomposition $(A_i,\phi_i)$ has 
the property $\mathcal{P}$ if for any $i \in \mathbb{N}$, element $f$ in ${A}_{i}$, 
$\epsilon > 0$, and $x \in X_i$ there exist $j \ge i$ and unitaries $u_l \in A_{j,l}$,
$l \in \{1,\ldots,k_j\}$ such that 
\begin{displaymath}
\left| \left| u_l \phi_{i,j}^{t,l}(f) u_l^* -  \left( \begin{array}{cc}
      f(x) & 0 \\
      0        & b_l
       \end{array} \right) \right| \right| < \epsilon
\end{displaymath}
for some appropriately sized $b_l$.  (Note that
$\mathrm{diag}(f(x_{0}),b_{l}) \in A_{{j,l}}$.)
\end{defn}

\vspace{2mm}
\noindent
We will prove in the sequel that $A$ as in Definition \ref{defn1} is simple if and only if
it has property $\mathcal{P}$.

\subsection{A characterisation of simplicity}
Proposition 2.1 of \cite{dnnp} gives some necessary and sufficent conditions 
for the simplicity of an AH algebra.  We will have occasion to apply these
in the proof of our main result, and so restate the said proposition 
in the particular case of an AH system with diagonal injective maps.
\begin{prop}\label{prop1.1}
Let $(A_i,\phi_i)$ be an AH system with diagonal injective maps, and set
$A = \lim_{i \to \infty}(A_i,\phi_i)$.  The following conditions are equivalent:
 \begin{itemize}
 \item[(i)] $A$ is simple;
 \item[(ii)] For any positive integer $i$ and any non-empty open subset $U$ of $X_{i}$, 
there is a $j_{0} \ge i$ such that for every $j \ge j_{0}$ and $l \in \{1,\ldots,k_j\}$
we have 
\[
(\mathrm{ep}_{ij}^l)^{-1}(U) = X_{j,l},
\]
where $(\mathrm{ep}_{ij}^l)^{-1}(U)$ denotes the union of the sets $\lambda^{-1}(U)$,
$\lambda \in \mathrm{ep}_{ij}^l$;
 \item[(iii)] For any non-zero element $a$ in ${A}_{i}$, there is a $j_{0} \ge i$ such that 
for every $j \ge j_{0}$, ${\phi}_{ij}(a)(x)$ is not zero, for every 
 $x$ in $X_{j}$.
 \end{itemize}
\end{prop}

\subsection{Paths between permutation matrices}\label{specialpath}

Given any permutation $\pi \in S_{n}$, let $U[\pi]$ denote the permutation matrix in $\mathrm{M}_n$ corresponding 
to $\pi$, that is, $U[\pi]$ is obtained from the identity of $\mathrm{M}_n$ by moving the $i^{\textrm{th}}$ row to 
the ${\pi (i)}^{\textrm{th}}$ row, for $i \in \{1,2, \ldots, n\}$. 
Any two permutation matrices are homotopic inside the unitary group $\mathcal{U}(\mathrm{M}_n)$ of $\mathrm{M}_n$,
but we want to define some particular homotopies for use in the sequel.
Let $\pi$ and $\sigma$ be elements of $S_{n}$, viewed as permutations of the canonical basis vectors
$e_1,\ldots,e_n$ of $\mathbb{C}^n$.  Let $R = \{e_{W,1},\ldots,e_{W,\mathrm{dim}(W)}\}$ be the set of
basis vectors upon which $\pi$ and $\sigma$ agree, and choose $\gamma \in S_n$ be such that 
\[
\gamma \sigma (v) = \gamma \pi (v) = v, \ \forall v \in R.
\]
Then $U[\gamma] U[\sigma]$ and $U[\gamma] U[\pi]$ fix $e_1,\ldots,e_{|R|}$.  
Put $W = \mathrm{span} \{e_1,\ldots,e_{|R|}\}$, and let $V$ be the orthogonal
complement of $W$.  There are a canonical unital embedding of $\mathrm{M}_{\mathrm{dim}(W)} \oplus
\mathrm{M}_{\mathrm{dim}(V)}$ into $\mathrm{M}_n$ and unitaries $u,v \in 
\mathcal{U}(\mathrm{M}_{\mathrm{dim}(V)})$ such that
\[
U[\gamma]U[\pi] = \mathbf{1}_{\mathrm{M}_{\mathrm{dim}(W)}} \oplus v; \ \ \ 
U[\gamma]U[\sigma] = \mathbf{1}_{\mathrm{M}_{\mathrm{dim}(W)}} \oplus u.
\]
Choose a homotopy $g(t)$ between $u$ and $v$ inside $\mathcal{U}(\mathrm{M}_{\mathrm{dim}(V)})$ 
--- $g(0) = v$ and $g(1) = u$ --- and put 
\[
u(t) = U[\gamma]^{-1}(\mathbf{1}_{\mathrm{M}_{\mathrm{dim}(V)}} \oplus g(t)).  
\]
Then $u(t)$ is a homotopy of unitaries between $U[\pi]$ and $U[\sigma]$ such that 
\[
u(t)(v) = U[\pi](v) = U[\sigma](v), \ \forall v \in R, \ \forall t \in [0,1].
\]

\subsection{Applications of Urysohn's Lemma}

\begin{lemma}\label{lem3.1}
Let $\sigma \in S_{n}$ be given.  There is a homotopy $u:[0,1] \to \mathcal{U}(M_{n})$ between $u(0) = 1_n$ 
and $u(1) = U[\sigma]$ which moreover has the following property:  for any complex numbers $\la_{1}, \la_{2}, \ldots, 
\la_{n}$ such that $\la_{i}$ = $\la_{\sigma (i)}$ for every $i=1, 2, \ldots, n$, we have
\begin{displaymath}
u(t) \left( \begin{array}{cccc}
          {\la}_{1} & 0   & \ldots & 0\\
         0     &  {\la}_{2} & \ldots & 0\\
         \vdots  &  \vdots & \ddots   &      \vdots\\
         0 &   0 &  \ldots &           {\la}_{n}
    \end{array} \right) u^*(t)
 =  \left( \begin{array}{cccc}
          {\la}_{1} & 0   & \ldots & 0\\
         0     &  {\la}_{2} & \ldots & 0\\
         \vdots  &  \vdots & \ddots   &      \vdots\\
         0 &   0 &  \ldots &             {\la}_{n}
    \end{array} \right), \ \forall t \in [0,1].
\end{displaymath}
\end{lemma}
\begin{proof}
Let us first consider the case that $\sigma$ is a $k$-cycle.  The hypotheses of the lemma guarantee that
the desired conclusion holds already for $t \in  \{0,1\}$.  Choose the homotopy $u(t)$ as in Subsection 
\ref{specialpath} by using our given value of $\sigma$ and setting $\pi$ equal to the identity element of
$S_n$.  The hypothesis $\lambda_i = \lambda_{\sigma(i)}$, $i \in \{1,\ldots,n\}$, implies that there is
a $\lambda \in \mathbb{C}$ such that for each $i \in \{1,\ldots,n\}$ which is not fixed by $\sigma$ 
we have $\lambda_i = \lambda$.  In other words, if one decomposes $\mathrm{diag}(\lambda_1,\ldots,\lambda_n)$
into a direct sum of two diagonal matrices using the decomposition $\mathbb{C}^n = W \oplus V$ --- $V$ and $W$ as in 
Subsection \ref{specialpath} --- then the direct summand corresponding to $V$ is scalar $k \times k$ matrix.  
By construction, $u(t) = v(t) \oplus 1_{n-k}$, with $v(t) \in \mathcal{U}(M_k)$.  It follows that $u(t)$
commutes with $\mathrm{diag}(\lambda_1,\ldots,\lambda_n)$ for each $t \in (0,1)$.

Now suppose now that $\sigma$ is any permutation on $n$ letters, and write $\sigma$ as a product 
of disjoint cycles:  $\sigma$ = $\sigma_{1} \sigma_{2} \ldots \sigma_{l}$. For each $j \in \{1,\ldots,l\}$, 
let $u_{i}(j)$ denote the unitary path between $U[\mathrm{id}]$ and $U[\sigma_{j}]$, constructed as in 
Subsection \ref{specialpath}.  Now 
\[
u(t) := u_1(t) u_2(t) \cdots u_l(t)
\]
is a path of unitaries with $u(0) = U[\mathrm{id}]$ and $u(1) = U[\sigma]$, and $u(t)$ commutes
with $\mathrm{diag}(\lambda_1,\ldots,\lambda_n)$ since each $u_j(t)$ does.
\end{proof}

\begin{lemma}\label{lemforthrm}
Let $\sigma$ be any permutation in $S_{n}$, and $A, B$ disjoint nonempty closed subsets of a metric space $X$. 
Let $\lambda_1,\ldots,\lambda_n:X \to \mathbb{C}$ be continuous.  Then, there exists a unitary $v \in M_{n}(\cx)$ such that
\begin{enumerate}
\item[(i)] $v(x) = 1_n$, $\forall x \in A$,
\item[(ii)] $v(x) = U[\sigma]$, $\forall x \in B$, and
\item[(iii)] $v(x)$ commutes with $\mathrm{diag}(\lambda_1(x),\ldots,\lambda_n(x))$ whenever $\lambda_i(x) = \lambda_{\sigma(i)}(x)$
for each $i \in \{1,\ldots,n\}$.
\end{enumerate}
\end{lemma}

\begin{proof}
Find a unitary path $u(t)$ connecting $U[\mathrm{id}]$ to $U[\sigma]$ using Lemma \ref{lem3.1},
so that $u(t)$ commutes with $\mathrm{diag}(\lambda_1(x),\ldots,\lambda_n(x))$ for each $t \in (0,1)$
and each $x \in X$ for which $\lambda_i(x) = \lambda_{\sigma(i)}(x)$, $i \in \{1,\ldots,n\}$. 
By Urysohn's Lemma, there is a continuous map $f:X \to [0,1]$ which is equal to zero on $A$ and 
equal to one on $B$.  It is straightforward to check that $v(x) := u(f(x))$ satisfies the conclusion
of the lemma.
\end{proof}

\begin{lemma}\label{tietze}
Let $Y$ be a closed subset of a normal space $X$, and let $f:Y \to \mathrm{S}^n$ be 
continuous.  Then, there is a continuous map 
$\tilde{f}$ from $X$ to the $n+1$ disk $D^{n+1}$ which extends $f$.
\end{lemma}

\begin{proof}
View $D^{n+1}$ as the cube $[0,1]^{n+1}$, with $S^{n}$ as its boundary.  Then
\[
f(y) = (f_1(y),\ldots,f_{n+1}(y)),
\]
where each $f_i:Y \to [0,1]$ is continuous.  Extend each $f_i$ to a continuous map
$\tilde{f}_i:X \to [0,1]$, and put
\[
\tilde{f}(y) = (\tilde{f}_1(y),\ldots,\tilde{f}_{n+1}(y)).
\]
\end{proof}

\section{The Main Theorem}

Let $A$ be a diagonal AH algebra with decomposition $(A_i,\phi_i)$, and
assume that the $\phi_i$ are injective.
In this section we will prove that if $A$ is simple then it has the 
property $\mathcal{P}$ (cf. Section 2).  (The converse also holds, but
is easier by far.)  Let us begin with an outline of our strategy, before
plunging headlong into the proof.

Assume first that $A$ is simple and $A_i = \mathrm{M}_{n_i}(\mathrm{C}(X_i))$,
so that there are no partial maps to contend with.
Let there be given a natural number $i$, an element $f$ of $A_i$,
a point $x_0 \in X_i$, and some $\epsilon > 0$.
Put $U = B_{\epsilon}(x_0)$.  By Proposition \ref{prop1.1} there is a $j_{0}$ 
with $j_{0} \ge i$ such that for any $j \ge j_{0}$, \[ X_{j} = {\la}_{1}^{-1}(U) \cup 
{\la}_{2}^{-1}(U) \cup \ldots \cup {\la}_{n}^{-1}(U), \] where 
\begin{displaymath}
\phi_{i,j}(f) =  \left( \begin{array}{cccc}
         f \circ \la_{1} &  0 &\ldots & 0  \\
       \vdots &   \vdots  &\ddots &  \vdots \\
        0 & 0     &\ldots &      f \circ \la_{n}
       \end{array} \right).
\end{displaymath}
On each closed subset $\la_{t}^{-1}(\overline{U})$, the range of the eigenvalue map $\la_{t}$ 
is within $\epsilon$ of $x_{0}$. To show that $A$ has the property $\mathcal{P}$, we require
a unitary $u$ in $A_{j}$ and an element $b_f \in \mathrm{M}_{n_j-n_i}(\mathrm{C}(X_j))$ such that 
\begin{displaymath}
\left| \left| u \phi_{i,j} (f) u^{*} - \left(
    \begin{array}{cc}
        f(x_{0}) & 0 \\
         0 & b_{f}
     \end{array}
         \right)  \right| \right| < \epsilon.
\end{displaymath}
We would like $u(y)$ to exchange the first and $t^{\mathrm{th}}$ diagonal entries of $\phi_{i,j}(f)$ 
whenever $y \in \lambda_t^{-1}(\overline{U})$, but this operation is unlikely to
be well-defined --- the sets $\lambda_1^{-1}(\overline{U}),\ldots,\lambda_{n}^{-1}(\overline{U})$ need not
be mutually disjoint.  The remainder of this section is devoted to overcoming this complication. 

Given positive integers $m \le n$, let us denote by $\mathrm{M}_{m}(\cy) \oplus 1_{n-m}$ the set of all $n \times n$ 
matrices of the form
\[ \left( \begin{array}{cc}
        a & 0 \\
         0 & 1
       \end{array}
         \right),   \]
where $a \in \mathrm{M}_{m}(\cy)$ and $1_{n-m} \in \mathrm{M}_{n-m}$.

\begin{thrm}\label{thrm}
Let there be given a diagonal $*$-homomorphism $\phi:\cx \to \mathrm{M}_n(\cy)$ with the eigenvalue pattern 
$\{{\la}_{1}, {\la}_{2}, \ldots, {\la}_{n}\}$, a point $x_{0}$ in $X$, an element $f$ of $\cx$, and 
a tolerance $\e>0$. Choose $\eta>0$ such that $|f(x)-f(y)| < \epsilon$ whenever $d(x,y) <2 \eta$
($d$ is the metric on $X$).
Suppose that $F_{1}, \ldots, F_{m}$ are nonempty closed subsets of $Y$ ($m \le n$) 
such that $d(\la_{i}(y), x_{0}) < \eta$ whenever $y \in F_{i}$.

Then, there is a unitary $u$ in $\mathrm{M}_{m}(\cy) \oplus 1_{n-m}$ and an element $b \in 
\mathrm{M}_{n-m-1}(\cy)$ such that 
for each $y \in \cup_{i=1}^m F_i$ we have
\begin{equation}\label{approximation}
\left| \left| u(y) \phi (f)(y) u^{*}(y) - \left(
    \begin{array}{ccccc}
        f(x_{0}) & 0 & 0 & \ldots & 0\\
         0 & b(y) & 0 & \ldots & 0 \\
         0 & 0 & \la_{m+1}(y) & \ldots & 0 \\
         \vdots & \vdots & \ddots & \vdots & \vdots \\
         0 &    0 &  0 & \ldots &            \la_{n}(y)
       \end{array}
         \right)  \right| \right| <2 \epsilon.
\end{equation}
\end{thrm}
(Note that if $n-m-1$ = $0$, then there is no $b$ in the matrix in (\ref{approximation}).) 

\begin{proof}
Let $\rho$ denote the metric on $Y$.
Choose $\delta>0$ such that $d(\la_{i}(x), \la_{i}(y))< \eta$ whenever $\rho (x,y) \le \delta$, $i \in \{1,\ldots,n\}$.
For each $1 \le i \le m$, 
 \[ |f \circ \la_{i}(y) - f(x_{0})| < \e, \ \textrm{ for all } y \in \overline{B_{\delta}(F_{i})}. \]
Set $\varepsilon_{i}(y)$ = $f \circ \la_{i}(y) - f(x_{0})$ for all $y $ in $\overline{B_{\delta}(F_{i})}$. 
Then, $\varepsilon_{i}$ is a continuous map from $\overline{B_{\delta}(F_{i})}$ to the disk of radius $\e$ 
in the complex plane. By Lemma \ref{tietze}, $\varepsilon_{i}$ can be extended to a continuous function 
from $Y$ to the complex plane such that $||\varepsilon_{i}||$ $\le $ $\e$ (let us also denote this 
extension map by $\varepsilon_{i}$). For $m < i\le n$, set $\varepsilon_{i}$ = $0$. 

For each $i \in \{1,\ldots,n\}$, put $g_i = f \circ \la_i - \epsilon_i$, so that $g_i \in \cx$.  Set
\[ 
g = \mathrm{diag}(g_{1}, g_{2}, \ldots, g_{n}).
\]
Then, for each $i \in \{1,\ldots,m\}$ and $y \in \overline{B_{\delta}(F_{i})}$, we have
\[ 
g_{i}(y) = f(x_{0}); 
\]
if $i \in \{m+1,\ldots,n\}$, then $g_i = f \circ \la_i$.
For any unitary $u \in \mathrm{M}_n(\cy)$ we have 
\[
|| u \phi(f) u^{*} - u g u^{*}|| = ||\mathrm{diag}(\varepsilon_{1}, \ldots, \varepsilon_{n}) || <2 \e.
\]
We have therefore reduced our problem to proving the following claim:

\vspace{2mm}
\noindent
{\bf Claim:} There is a unitary $u$ in $\mathrm{M}_{m}(\cy) \oplus 1_{n-m}$ such that 
\[ ugu^{*} = \left(
    \begin{array}{ccccc}
        f(x_{0}) & 0 & 0 & \ldots & 0\\
         0 & b(x) & 0 & \ldots & 0 \\
         0 & 0 & g_{m+1}(x) & \ldots & 0 \\
         \vdots & \vdots & \vdots & \ddots & \vdots \\
         0 &    0 &  0 & \ldots &            g_{n}(x)
       \end{array}
         \right), \ \forall x \in \bigcup_{i=1}^{m}F_{i},   \]
where $b \in \mathrm{M}_{n-m-1}(\cy)$. 

\vspace{2mm}
\noindent
{\bf Proof of claim:}
We will assume that for some $1 \leq k < m$ there is a unitary $u_k \in 
\mathrm{M}_k(\cy) \oplus 1_{n-k}$ such that
\begin{equation}\label{claim} 
u_kgu_k^{*} = \left(
    \begin{array}{ccccc}
        f(x_{0}) & 0 & 0 & \ldots & 0\\
         0 & b(x) & 0 & \ldots & 0 \\
         0 & 0 & g_{k+1}(x) & \ldots & 0 \\
         \vdots & \vdots & \vdots & \ddots & \vdots \\
         0 &    0 &  0 & \ldots &            g_{n}(x)
       \end{array}
         \right), \ \forall x \in \bigcup_{i=1}^{k}F_{i},   
\end{equation}
and then prove that the same statement holds with $k$ replaced by $k+1$.
Since (\ref{claim}) clearly holds when $k=1$ --- just take $u_1$ to be 
the identity matrix of $\mathrm{M}_{n}(\cy)$ --- this recursive argument will prove
our claim.

Assume that (\ref{claim}) holds for some $k < m$.  
Put $B = F_{k+1}$ and $A = Y \backslash B_{\delta}(B)$. Apply Lemma \ref{lemforthrm}
with these choices of $A$ and $B$ and with $\sigma = (1 \ k+1)$ to obtain a unitary
$v \in \mathrm{M}_n(\cy)$.  We then have that $v=1_n$ on $A$ and $v=U[(1 \  k+1)]$ on $B$.
Inspecting the construction of $v$, we find that it has the following form: 
\begin{equation}\label{vfactor} 
v(y) = U[(2 \  k+1)] \left( \begin{array}{cc} v^{'}(y)& 0 \\ 0 & 1 \end{array}
\right)U[(2 \ k+1)], \forall y \in Y, 
\end{equation}
where $v^{'}(y)$ is a unitary matrix in $\mathrm{M}_{2}(\cy)$ equal to $1_2$ on $A$ and equal
to $U[(1  2)]$ on $B$.  Define $u_{k+1} := vu_k$.  

Let us show that $u_{k+1}$ satisfies the requirements of the claim.
It is clear that $u_{k+1}$ is an element of $\mathrm{M}_{k+1}(\cy) \oplus 1_{n-k-1}$.
First suppose that $y \in B$, so that $g_{k+1}(y) = f(x_0)$.  Since $u_{k+1} = v u_k$ we have
\begin{eqnarray*} u_{k+1}(y)g(y)u_{k+1}^*(y) & = & U[(1 \ k+1)] 
\left( \begin{array}{cccc}
      c(y)& 0 & \ldots & 0 \\
      0 & f(x_0) & \ldots & 0\\
      \vdots & \vdots & \ddots & \vdots \\
      0 & 0   & \ldots &   g_{n}(y)
      \end{array} \right) U[(1 \ k+1)] \\
& = &  \left( \begin{array}{ccccc}
         f(x_{0}) &  0 & 0 & \ldots & 0  \\
        0  &  b(y)  & 0 &  \ldots & 0 \\
        0 & 0  &  g_{k+2}(y) & \ldots & 0 \\
        \vdots &   \vdots  & \vdots & \ddots &  \vdots \\
        0 & 0   & 0  &\ldots &      g_{n}(y)
       \end{array} \right), 
\end{eqnarray*} 
for some  $c(y),b(y) \in \mathrm{M}_{k}$.

Now suppose that $y \in \bigcup_{i=1}^{m}F_{i} \backslash B_{\delta}(B) \subseteq Y \backslash B_{\delta}(B)$.
In this case $v(y)=1_n$ and $u_{k+1}(y) = u_k(y)$ and there is nothing to prove.

Finally, suppose that $y \in (B_{\delta}(B) \backslash B) \cap (\bigcup_{i=1}^{m}F_{i})$.  As in the case $y \in B$,
we have $g_{k+1}(y) = f(x_0)$.  From (\ref{claim}) and this last fact we have 
\[ U[(2 \ k+1)]u_k(y)g(y)u_k^{*}(y)U[(2 \ k+1)] =
   \left( \begin{array}{cccccc}
         f(x_{0}) &  0 &    0 & 0 & \ldots & 0  \\
          0 & f(x_{0}) & 0  & 0 & \ldots &0 \\
          0 & 0 &  d(y) & 0 & \ldots & 0\\  
          0 & 0 & 0 & g_{k+2} (y) & \ldots & 0\\
         \vdots &   \vdots  & \vdots & \vdots & \ddots &  \vdots \\
        0 & 0 & 0   & 0 & \ldots &      g_{n}(y)
       \end{array} \right)
\]
for some $d(y) \in M_{k-1}$.  Since the upper left $2 \times 2$ corner of
the matrix above is scalar, the entire matrix commutes with $v^{'}(y) \oplus 1_{n-2}$.
It follows that $u_{k+1}(y) g(y) u_{k+1}^{*}(y)$ is equal to
\begin{displaymath}
   U[(2 \  k+1)]
      \left( \begin{array}{cccccc}
         f(x_{0}) &  0 &    0 & 0 & \ldots & 0  \\
          0 & f(x_{0}) & 0  & 0 & \ldots &0 \\
          0 & 0 &  b_{1}(y) & 0 & \ldots & 0\\  
          0 & 0 & 0 & g_{k+2} (y) & \ldots & 0\\
         \vdots &   \vdots  & \vdots & \vdots & \ddots &  \vdots \\
        0 & 0 & 0   & 0 & \ldots &      g_{n}(y)
       \end{array} \right)
  U[(2 \  k+1)].    
\end{displaymath}
Computing this product yields a matrix of the form
\begin{displaymath}
\left( \begin{array}{ccccc}
         f(x_{0}) &  0 & 0& \ldots & 0   \\
        0  &  b(y)  & 0 &  \ldots & 0 \\
        0 & 0  &  g_{k+2}(y) & \ldots & 0 \\
        \vdots &   \vdots  & \vdots & \ddots &  \vdots \\
        0 & 0   & 0  &\ldots &      g_{n}(y)
       \end{array} \right),
\end{displaymath} 
as required.
\end{proof}

For any C*-algebra $B$ there is an isomorphism
\[
\pi:B \otimes \mathrm{M}_n \to \mathrm{M}_n(B)
\]
given by
\begin{equation}\label{tensorform}
\pi(b \otimes (a_{ij})) = \left( \begin{array}{ccc} b a_{11} & \cdots & b a_{1n} \\
\vdots & \ddots & \vdots \\ b a_{n1} & \cdots & b a_{nn} \end{array} \right).
\end{equation}


\begin{prop}\label{matrixtensor} 
Suppose that $\phi:\mathrm{M}_m(\cx) \to \mathrm{M}_{nm}(\cy)$ is a diagonal $*$-homomorphism with 
$\textrm{ep}(\phi)$ = $\{ \la_{1}, \la_{2}, \ldots, \la_{n}\}$.  Let $\tilde{\phi}:\cx \to \mathrm{M}_n(\cy)$ 
be the diagonal $*$-homomorphism given by 
\[ \tilde{\phi} (f) = \left( \begin{array}{cccc}
      f \circ \la_{1} & 0  &  \ldots & 0   \\
      0 &  f \circ \la_{2} & \ldots & 0  \\
      \vdots & \vdots &  \ddots &  \vdots \\
      0  &   0 & \ldots &      f \circ \la_{n}
       \end{array} \right).  
\]
Then, $\tilde{\phi} \otimes \mathbf{id}_{\mathrm{M}_m}$: $\cx \otimes \mathrm{M}_m \ \longrightarrow \mathrm{M}_n(\cy) 
\otimes \mathrm{M}_m$ is unitarily equivalent to $\phi$.  
\end{prop}

\begin{proof}

On the one hand we have
\begin{displaymath}
\tilde{\phi} \otimes \mathbf{id}_{\mathrm{M}_m}(f \otimes (c_{ij})) = \left( \begin{array}{cccc}
      f \circ \la_{1} & 0  &  \ldots & 0   \\
      0 &  f \circ \la_{2} & \ldots & 0  \\
      \vdots & \vdots &  \ddots &  \vdots \\
      0  &   0 & \ldots &      f \circ \la_{n}
       \end{array} \right) \otimes (c_{ij}), 
\end{displaymath}
while on the other we have
\[
\phi(f \otimes (c_{ij})) = \left( \begin{array}{cccc}
      (f \circ \la_{1}) \otimes (c_{ij}) & 0  &  \ldots & 0   \\
      0 &  (f \circ \la_{2}) \otimes (c_{ij}) & \ldots & 0  \\
      \vdots & \vdots &  \ddots &  \vdots \\
      0  &   0 & \ldots &      (f \circ \la_{n}) \otimes (c_{ij})
       \end{array} \right).
\] 
With the identifications $\cx\otimes \mathrm{M}_n \cong \mathrm{M}_n(\cx)$ and $\mathrm{M}_n(\cy) 
\otimes \mathrm{M}_m \cong \mathrm{M}_{nm}(\cy)$ given by (\ref{tensorform}) in mind, one sees
that 
\begin{displaymath}
\mathrm{Ad}(U[\pi]) \circ (\tilde{\phi} \otimes \mathbf{id}_{\mathrm{M}_m}) = \phi, 
\end{displaymath}
where $\pi$ is the permutation in $S_{nm}$ and given by $\pi(kn + i )$ = $(i-1)m + k+1$ for 
$k = 0,1,2, \ldots, m-1$ and $i$ = $1,2, \ldots, n$. 
\end{proof}

\begin{cor}\label{corclarify}
Let $\phi:\mathrm{M}_m(\cx) \to \mathrm{M}_{nm}(\cy)$ be a diagonal $*$-homomorphism 
with eigenvalue pattern $\{{\la}_{1}, {\la}_{2}, \ldots, {\la}_{n}\}$, and let $\epsilon>0$
be given.  Let $f$ be any element of $\mathrm{M}_m(\cx)$ and choose $\eta >0$ such that 
\begin{displaymath}
 ||f(x)-f(y)|| < \epsilon  \ \ \mathrm{ whenever } \ \ d(x,y) <2 \eta
\end{displaymath}
($d$ is the metric on $X$). Let $U$ be the open ball centered at $x_{0} \in X$ with radius $\eta$, and
suppose that 
\begin{displaymath}
Y = {\la}_{1}^{-1}(U) \cup {\la}_{2}^{-1}(U) \cup \ldots \cup {\la}_{n}^{-1}(U).
\end{displaymath}
Then there is a unitary $u$ in $\mathrm{M}_{nm}(\cy)$  and element $b \in 
\mathrm{M}_{nm-m}(\cy)$ such that 
\begin{displaymath}
\left| \left| u \phi (f) u^{*} - \left(
    \begin{array}{cc}
        f(x_{0}) & 0 \\
         0 & b
       \end{array}
         \right)  \right| \right| < \epsilon.
\end{displaymath}
\end{cor}

\begin{proof}
By Proposition \ref{matrixtensor}, we can assume that $m$ = $1$. 
Let $F_{i}$ be the closure of $ {\la}_{i}^{-1}(U)$ for each $i$. 
Now, apply Theorem \ref{thrm}.  Since $Y$ = $\bigcup_{i=1}^{n} F_{i}$,
we are done.
\end{proof}

Now, we are ready to prove the main theorem of this section.

\begin{thrm} \label{main}
Let $A = \underrightarrow{\lim} (A_{i}, \phi_{i})$ be a unital diagonal AH algebra. 
Then, $A$ is simple if and only if $A$ has the property $\mathcal{P}$ of Definition \ref{defn1}.
\end{thrm}

\begin{proof}
Suppose that $A$ has property $\mathcal{P}$.  Let $f \in A_{i} $ be nonzero, so that there 
is a point $x_{0}$ in $X_{i}$ such that $f(x_{0}) \ne 0$. By the definition of property 
$\mathcal{P}$, there are an integer $j > i$ and unitaries $u_l \in A_{j,l}$, $l \in \{1,\ldots,k_j\}$
such that
\begin{displaymath}
\left| \left| u_l \phi_{i,j}^{t,l}(f) u_l^* -  \left( \begin{array}{cc}
      f(x_0) & 0 \\
      0        & b_l
       \end{array} \right) \right| \right| < \epsilon
\end{displaymath}
for some appropriately sized $b_l$.
We may assume that $\epsilon < ||f(x_0)||$, so that ${\phi}_{ij}(f)$ is nowhere zero. 
This implies that the ideal of $A_j$ generated by $\phi_{ij}(f)$ is all of $A_j$, and
that the ideal of $A$ generated by the image of $f$ is all of $A$.  Since $f$ was arbitrary,
$A$ is simple.

Now assume that $A$ is simple, and let $f \in A_{i,t}$ be nonzero. 
Recall that the $\phi_i$ may be taken to be injective.
By Proposition \ref{prop1.1} there exists, for each $x_0 \in X_{i,t}$ and $\epsilon>0$,
a $j_0 > i$ with the following property:  for every $j \geq j_0$ and every $l \in \{1,\ldots,k_j\}$
we have
\[
\textrm{ep}^{-1}(\phi_{i,j}^{t,l})(B_{\delta}(x_{0})) = X_{j,l},
\]
where $\delta$ is some positive number such that 
\[
d(f(x),f(y)) < \epsilon \ \ \mathrm{ whenever } \ \ d(x,y) < 2 \delta.
\]
As pointed out in Section 2, the map $\phi_{i,j}^{t,l}$ may be viewed as
a diagonal map from $A_{i,t}$ into the cut-down of $A_{j,l}$ by the projection
$\phi_{i,j}^{t,l}(1)$;  any unitary $u$ in this corner of $A_{j,l}$ gives
rise to a unitary $\tilde{u}$ in $A_{j,l}$ by setting $\tilde{u}= u \oplus (1_{A_{j,l}} - 
\phi_{i,j}^{t,l}(1)$.  Combining this observation with Corollary \ref{corclarify} 
we see that there exists, for each $l \in \{1,\ldots,k_j\}$, a unitary $u_l \in A_{j,l}$
such that 
\begin{displaymath}
\left| \left| u_l \phi_{i,j}^{t,l}(f) u_l^* -  \left( \begin{array}{cc}
      f(x_0) & 0 \\
      0        & b_l
       \end{array} \right) \right| \right| < \epsilon.
\end{displaymath}
Thus, $A$ has property $\mathcal{P}$, as desired.
\end{proof}


\section{Stable Rank}



\begin{thrm}\label{stable}
Let $A = \underrightarrow{\lim} (A_{i}, \phi_{i})$ be a simple unital diagonal 
AH algebra.  Then, $A$ has stable rank one.
\end{thrm}

Before proving Theorem \ref{stable}, we situate it relative to
other results on the stable rank of general approximately homogeneous (AH) algebras.
Recall that an AH algebra is an inductive limit $C^*$-algebra
$A= \lim_{i \to \infty}(A_i,\phi_i)$, where
\begin{equation}\label{decomp}
A_i = \bigoplus_{l=1}^{n_i} p_{i,l}(\mathrm{C}(X_{i,l}) \otimes \mathcal{K})p_{i,l}
\end{equation}
for compact connected Hausdorff spaces $X_{i,l}$, projections $p_{i,l} \in
\mathrm{C}(X_{i,l}) \otimes \mathcal{K}$, and natural numbers $n_i$.  If
$A$ is separable, then one may assume that the $X_{i,l}$ are finite CW-complexes
(\cite{Bl2}, \cite{Go1}).  The inductive system $(A_i,\phi_i)$ is referred to 
as a {\it decomposition} for $A$.  All AH algebras in this paper are assumed to 
be separable.

If an AH algebra $A$ admits a decompostion as in (\ref{decomp}) for which 
\[
\mathrm{max}_{1 \leq l \leq n_i} \ \left\{ \frac{\mathrm{dim}(X_{i,1})}{\mathrm{rank}(p_{i,1})}, \ldots,
\frac{\mathrm{dim}(X_{i,n_i})}{\mathrm{rank}(p_{i,n_i})} \right\} \stackrel{i \to \infty}{\longrightarrow} 0,
\]
then we say that $A$ has {\it slow dimension growth}.  
Theorem 1 of \cite{bdr} states that every simple unital AH algebra with
slow dimension growth has stable rank one.  Villadsen in \cite{jes1} constructed
simple diagonal AH algebras which do not have slow dimension growth, but 
which do have stable rank one;  the converse of \cite[Theorem 1]{bdr} does
not hold.  There are in fact a wealth of simple diagonal AH algebras without slow dimension
growth which exhibit all sorts of interesting behaviour (cf. \cite{To1}, \cite{To2}, \cite{To3}), whence Theorem
\ref{stable} is widely applicable. 

Simple AH algebras may have stable rank strictly greater than one, and there is 
reason to believe that Theorem \ref{stable} is quite close to being best possible.
One might be able to generalise our result to the setting of AH algebras where the 
projections $\phi_{i,j}(p_{i,l})$ appearing in (\ref{decomp}) can be decomposed
into a direct sum of a trivial projection $\theta_j$ and a second projection $q_j$
such that $\tau(q_j) \to 0$ as $j \to \infty$ for any trace $\tau$.  Otherwise, one finds
oneself in a situation very similar to the construction of \cite{jes2}, where the stable 
rank is always strictly greater than one.



Let us now prepare for the proof of Therem \ref{stable}.

\begin{lemma}\label{lem2.1}
Let $a \in \mathrm{M}_m(\cx)$ be block diagonal, i.e., 
\begin{displaymath}
a = \left( \begin{array}{cccc}
      0 & 0   & \ldots &   0 \\
      0 & a_{1} & \ldots & 0\\
      \vdots & \ddots & \vdots \\
      0 & 0  &  \ldots &  a_{n}
       \end{array} \right),
\end{displaymath}
where $a_{i} \in \mathrm{M}_{k_{i}}(\cx)$ for natural numbers $k_1,\ldots,k_n$.  If the size of 
the matrix $0$ in the upper left-hand corner of $a$ is strictly greater than $\max_{1 \le i \le n}k_{i}$, 
then $a$ can be approximated arbitrarily closely by invertible elements in $\mathrm{M}_m(\cx)$.
\end{lemma}

\begin{proof}
Let $\epsilon>0$ be given, and let $k$ denote the size of the matrix $0$ in the upper left-hand corner of $a$. 
Since $k > k_i$, $i \in \{1,\ldots,n\}$, there is a permutation matrix $U$ such that
\begin{displaymath}
Ua = \left( \begin{array}{cccc}
      0 & *   & \ldots &   * \\
      0 & 0 & \ldots & *\\
      \vdots & \ddots & \vdots \\
      0 & 0  &  \ldots &  0
       \end{array} \right),
\end{displaymath}
is nilpotent.  As was proved in \cite{ro}, every nilpotent element in a unital C$^*$-algebra can be approximated
arbitrarily closely by invertible elements.  We may thus find an invertible element $b
\in \mathrm{M}_{m}(\cx)$ such that $||Ua-b|| < \epsilon$, and 
\[
|| a - U^{-1}b || = || U^{-1}Ua - U^{-1}b || \leq ||U^{-1}|| \cdot || Ua -b || < \epsilon.
\]
The lemma now follows from the fact that $U^{-1}b$ is invertible.
\end{proof}

It easy to prove that $a \in \mathrm{M}_n(\cx)$ is invertible if and only if $a(x)$ is
invertible for each $x \in X$.  The proof of the next lemma is also straightforward.


\begin{lemma}\label{rem2.2}
Let $p,q$ be orthogonal projections in a C$^*$-algebra $A$, and let $\epsilon>0$ be given. 
If elements $a$ and $b$ in $A$ can be approximated to within $\epsilon$ by invertible elements 
in $p A p$ and $q A q$, respectively, then $a + b$ can be approximated to within $\epsilon$ by 
an invertible element in $(p+q)A(p+q)$.
\end{lemma}

Now, we are ready to prove Theorem \ref{stable}.
\begin{proof}[Proof of Theorem \ref{stable}]
Since every element in $A$ can be approximated arbitrarily closely by elements 
in $\bigcup_{i=1}^{\infty} A_{i}$, it will suffice to prove for any $\epsilon>0$
and any $a \in A_{i}$, there is an invertible element
in $A$ whose distance to $a$ is less than $\epsilon$.
(Note that we are using the injectivity of the $\phi_i$ to identify $A_{i}$ with its image in $A$.)

By Lemma \ref{rem2.2}, we may assume that $A_{i} = \mathrm{M}_{n_{i}}(\mathrm{C}(X_{i}))$.  We also assume  
that $a$ is not invertible. By the comment preceding Lemma \ref{rem2.2}, there is a point $x_{0} 
\in X_{i}$ such that det($a(x_{0})$)= $0$. There are permutation matrices $u,v
\in \mathrm{M}_{n_{i}}$ and a matrix $c \in \mathrm{M}_{n_{i}-1}$ such that
\begin{displaymath}
u a(x_{0}) v = \left( \begin{array}{cc}
      0 & 0  \\
      0 & c
      \end{array} \right),
\end{displaymath}
Let $b$ denote the element $uav$.  Following the lines of the proof of Lemma
\ref{lem2.1}, it 
will suffice to prove that $b$ can be approximated to within $\epsilon$ by an invertible element
of $A$. 

For each $j>i$, $\phi_{i,j}(b)$ is a tuple of $k_{j}$ elements. If each 
coordinate of $\phi_{i,j}(b)$ can be approximated to within $\epsilon$ by 
an invertible element in the corner of $A$ generated by the unit of
$A_{j,l}$, then $\phi_{i,j}(b)$ can be approximated to within $\epsilon$ 
by an invertible element of $A$.  We may therefore assume that $A_{j}=
\mathrm{M}_{n_{j}}(\mathrm{C}(Y_{j})$, and concern ourselves with
proving that $\phi_{i,j}(b)$ is approximated to within $\epsilon$ by an invertible element 
in $A$.

By Theorem \ref{main}, there exist an integer $j>i$, a unitary $w \in A_{j}$ 
and an element $b'$ such that 
\begin{displaymath}
 \left| \left| w \phi_{ij} (b)w^{*} - \left( \begin{array}{cc}
                             b(x_{0}) & 0 \\
                             0        & b' 
                             \end{array} \right) \right| \right| < \epsilon/2.
  \end{displaymath}
We have
\begin{displaymath}
  \left( \begin{array}{cc}
         b(x_{0}) & 0 \\
         0        & b' 
         \end{array} \right) = 
        \left( \begin{array}{cc}
          0       & 0 \\
          0       & b'' 
         \end{array} \right),
  \end{displaymath}
where $b'' = \mathrm{diag} (c, b')$.  Put $d = w \phi_{i,j} (b)w^{*}$, and note that 
it will suffice to
prove that $\phi_{j,m}(d)$ is approximated to within $\epsilon/2$ by an invertible element in $A_m$ for
some $m>j$.  

Since $A$ is simple, there is an integer $m>j$ large 
enough so that, for each $t \in \{1,\ldots,k_m\}$, either the number of the eigenvalue maps of 
$\phi_{j,m}^{1,t}$ counted with multiplicity is strictly larger than the size of the matrix $b''$, 
or the image of $\phi_{j,m}^{1,t}$ is finite-dimensional.  In the latter case, $\phi_{j,m}^{1,t}(d)$
is approximated to within $\epsilon$ by an invertible element in the image of $\phi_{j,m}^{1,t}$
since finite-dimensional C$^*$-algebras have stable rank one, so we may assume that the number of eigenvalue
maps of $\phi_{j,m}^{1,t}$, counted with multiplicity, is strictly larger than the size of the matrix $b''$.
Then, $\phi_{j,m}^{1,t}(d)$ is unitarily equivalent to
\begin{displaymath}
  \left( \begin{array}{cccc}
         0        & 0    & \ldots & 0\\
         0        & b_{1} & \ldots & 0 \\
         \vdots  & \vdots & \ddots & \vdots \\
        0        &  0      & \ldots & b_{l} 
         \end{array} \right)
\end{displaymath}
inside $A_{m,t}$, where $b_{k}$ is the composition of $b''$ and the $k^{\textrm{th}}$ eigenvalue map 
of $\phi_{j,m}^{1,t}$, and the size of the matrix $0$ at the upper left-hand corner of the above matrix 
is strictly bigger than the size of the matrix $b_{k}$ for every $k$. By Lemma \ref{lem2.1}, the matrix 
above can be approximated to within $\epsilon/2$ by an invertible in $A_{m,t}$, as required.
\end{proof}

\begin{cor}\label{taf}
Let $A$ be a simple unital diagonal AH algebra.  If $A$ has real rank zero 
and weakly unperforated $\mathrm{K}_0$-group, then $A$ is tracially AF.
\end{cor}

\begin{proof}
By Theorem \ref{stable}, $A$ has stable rank one.  The corollary then follows from
a result of Lin (\cite[Theorem 2.1]{lin}).
\end{proof}

Corollary \ref{taf} applies, for instance, to simple unital diagonal AH algebras for which the spaces $X_i$
in some diagonal decomposition for $A$ are all contractible.  This contractibility hypothesis may seem
strong, but it does not substantially restrict the complexity of $A$;  if one wants to classify
all such $A$ via $\mathrm{K}$-theory and traces, then the additional assumption of very slow 
dimension growth and the full force of \cite{EGL} and \cite{G} are required;  the collection of all
such cannot cannot be classified by topological $\mathrm{K}$-theory and traces alone (\cite{To2}).

\end{document}